\documentclass[twoside]{article}
\usepackage{graphics}
\usepackage{graphicx}
\usepackage{latexsym}
\usepackage{amssymb}
\usepackage{amsmath}
\usepackage{BSPMsample}
\numberwithin{equation}{section}
\begin{document}

\title[Characterization of Spherical and Plane Curves Using Rotation Minimizing Frames]{Characterization of Spherical and Plane Curves Using Rotation Minimizing Frames\footnote{This is a pre-print of article accepted for publication in 2021 in the Boletim da Sociedade Paranaense de Matem\'atica. The final authenticated version is available online at http://dx.doi.org/10.5269/bspm.49075. (See also the journal's page http://www.spm.uem.br/bspm/pdf/next/401.pdf)}}
%between [] goes the title that will appear on the top of every
%odd page, between {} goes the real title.
%\thanks{...} put your grant

\author[L. C. B. da Silva]{Luiz C. B. da Silva}
%between [] goes the authors that will appear on the top of every
%even page, between {} goes the real author.

\address{Luiz C. B. da Silva,\\Department of Physics of Complex Systems, \\
Weizmann Institute of Science, \\ Rehovot 7610001, Israel.}
\email{luiz.da-silva@weizmann.ac.il}

\maketitle

\begin{abstract}
In this work, we study plane and spherical curves in Euclidean and Lorentz-Minkowski 3-spaces by employing rotation minimizing (RM) frames. By conveniently writing the curvature and torsion for a curve on a sphere, we show how to find the angle between the principal normal and an RM vector field for spherical curves. Later, we characterize plane and spherical curves as curves whose position vector lies, up to a translation, on a moving plane spanned by their unit tangent and an RM vector field. Finally, as an application, we characterize Bertrand curves and slant helices as curves whose so-called natural mates are spherical and general helices, respectively.
\end{abstract}
%%put here Key words


\tableofcontents

%% put here the subject class
\2010mathclass{53A04, 53A55, 53B99}

\section{Introduction}

The usual way of studying curves is by means of the Frenet frame. But since its principal normal always points to the center of curvature, it may result in unnecessary rotation. In this respect, the consideration of \emph{rotation minimizing frames} (RM frames, for short) may be of special interest \cite{BishopMonthly,KlokCAGD1986}. The basic idea is that the normals should rotate only the necessary amount to remain normal to the tangent. The goal of this work is to study some properties of spherical and plane curves using RM frames. We show how to find the angle between the principal normal and an RM vector for a $C^3$ space curve (the derivative of this angle gives the torsion). This is done by using a convenient expression for the curvature function and torsion of a spherical curve. Subsequently, we  generalize these expressions for a generic curve by using  osculating spheres. In addition, we address the problem of characterizing curves whose position vector lies on a moving plane spanned by their unit tangent and an RM vector field and prove that they are precisely the plane and spherical curves. Finally, as an application of spherical curves and RM frames, we characterize Bertrand curves as the curves whose so-called natural mates \cite{DeshmukhTJM2018} are spherical. A notion of mate curves based on RM frames is also introduced and allows us to reverse the above mentioned association between spherical and Bertrand curves.

The remaining of this work is divided as follows. In Section 2, we introduce RM frames and some geometric background. In Section 3, we describe the behavior a twisted curve near its osculating sphere. In Section 4, we establish a characterization of plane and spherical curves using RM frames. In Section 5, we present {new characterizations of Bertrand curves and helices. Finally, in Section 6, we establish similar results for curves in the Lorentz-Minkowski space.}

\section{Preliminaries}
 
Denote by $\mathbb{E}^3$ the 3d Euclidean space, i.e., $\mathbb{R}^3$ equipped with the standard metric $\langle x,y\rangle=\sum_{i=1}^{\,3}x_iy_i$. Given a regular curve $\alpha:I\rightarrow \mathbb{E}^3$ parameterized by its arc-length $s$, $\Vert\alpha'(s)\Vert=1$, we may equip it with the  Frenet frame $\{\mathbf{t},\mathbf{n},\mathbf{b}\}$, whose equations of motion are $\mathbf{t}'=\kappa\,\mathbf{n}$, $\mathbf{n}'=-\kappa\,\mathbf{t}+\tau\,\mathbf{b}$, and $\mathbf{b}'=-\tau\,\mathbf{n}$, where $\kappa$ and $\tau$ are the curvature and torsion, respectively \cite{Kreyszig1991,Kuhnel2010}. If $\kappa>0$ and $\tau\not=0$, we say $\alpha$ is a \emph{twisted curve}. If $\kappa(s^*)=0$, $\alpha(s^*)$ is called an \emph{inflection point}.

Alternatively, one may consider an orthonormal frame $\{\mathbf{t},\mathbf{n}_1,\mathbf{n}_2\}$ with the additional property of $\mathbf{n}_i$ rotating around the unit tangent $\mathbf{t}$ only. The equation of motion of a {\it rotation minimizing} (RM) moving frame is 
\begin{equation}
\mathbf{t}'(s)=\kappa_{1}(s)\,\mathbf{n}_1(s)+ \kappa_{2}(s)\,\mathbf{n}_2(s),\,
\mathbf{n}_i'(s)=-\kappa_{i}(s)\,\mathbf{t}(s)\\
.\label{eq::BishopEqs}
\end{equation}
Now, writing $\mathbf{n}_1=\cos\theta\,\mathbf{n}-\sin\theta\,\mathbf{b}$ and $\mathbf{n}_2=\sin\theta\,\mathbf{n}+\cos\theta\,\mathbf{b}$, the curvatures $\kappa_1,\,\kappa_2$ relate to $\kappa$ and $\tau$ by 
$
\kappa_1(s)+\mathrm{i}\kappa_2(s) = \kappa(s)\,\mathrm{e}^{\mathrm{i}\theta(s)}$ and $\theta'(s)= \tau(s)$ \cite{BishopMonthly,GuggenheimerCAGD1989}. An advantage of RM  frames is that they can be globally defined even if $\kappa=0$ at some points \cite{BishopMonthly}.  On the other hand, they are not unique, since any rotation of $\mathbf{n}_i$ on the normal plane still gives an RM field, i.e., $\theta$ is defined up to a constant. Nonetheless,  $\kappa_1$ and $\kappa_2$  determine $\alpha$ up to rigid motions \cite{BishopMonthly}. Finally, RM frames allow for a simple characterization of spherical and plane curves\footnote{An attempt to extend these ideas in order to characterize curves that lie on a surface is described in Ref. \cite{daSilvaJG2017}.}
\begin{thm}[Bishop \cite{BishopMonthly}]\label{thr::charSphCurves}
A regular $C^2$ curve $\alpha:I\to \mathbb{E}^3$ lies on a sphere of radius $r$ if and only if its normal development curve $(\kappa_1(s),\kappa_2(s))$ lies on a line not passing through the origin where the distance from the origin is $r^{-1}$. In addition, normal development curves passing through the origin characterize plane curves.
\end{thm}

To find an RM frame we may compute the angle between  a normal vector field and the principal normal. In order to be RM, we should have $\theta'=\tau$ \cite{GuggenheimerCAGD1989}. The drawback here is the need of a Frenet frame globally defined, i.e., no inflection point should be allowed, otherwise $\tau$ may be not defined. On the other hand, the problem is easy to solve for  $\alpha\subset \mathbb{S}^2(p,r)$. The normalized position vector $\mathbf{N}=\frac{1}{r}(\alpha-p)$ is RM: ${\rm d}\mathbf{N}/{\rm d}s=\alpha'/r$. (The curvature  associated with $\mathbf{N}$ is then $\kappa_1=-\frac{1}{r}$.) This is an important step in the implementation of the {\it double reflection method} for computing approximations of RM frames \cite{WangACMTOG}. Another remarkable feature of RM frames along spherical curves is that they are path independent, i.e., if two spherical curves connect $q_1$ to $q_2$ in $\mathbb{S}^2(p,r)$ and their normals at $q_1$ coincide, then their normals at $q_2$ must also coincide \cite{WangACMTOG}. Notice that for the remaining RM vector field $\mathbf{n}_2=\alpha'\times\frac{(\alpha-p)}{r}$, one has $\kappa_2=-\frac{J}{r}$, where $J=\langle\alpha-p,\alpha'\times\alpha''\rangle$ is the \emph{spherical curvature} \cite{SabanRendLincei1958}. By noting that $\alpha'' = -\frac{1}{r}\,\frac{\alpha}{r}+\frac{J}{r}\,\frac{\alpha}{r}\times\alpha'$, one sees that $J$ is a multiple of the geodesic curvature of $\alpha\subset\mathbb{S}^2(p,r)$: $J = r \,\kappa_g$.
In addition, the curvature and torsion of a spherical curve are 
$\kappa = \frac{1}{r}\sqrt{1+J^2}$ and $\tau = \frac{J'}{1+J^2}$ \cite{Kuhnel2010,SabanRendLincei1958}.
This allows us to find the angle $\theta$ between an RM and the principal normal as $\theta(s) = \arctan\,J(s)+\mbox{constant}$.

\section{Behavior of a twisted curve near an osculating sphere}

We say that $\alpha$ and $\beta$ in $\mathbb{E}^3$ have a \emph{contact of order $k$} if $\alpha(s_0)=\beta(s_0^*)$ and all the higher order derivatives, up to $k$, also coincide: $\alpha^{(i)}(s_0)=\beta^{(i)}(s_0^*)$ for $1\leq i \leq k$. For example, the tangent line has a contact of order 1 with its reference curve, while the osculating circle has a contact of order 2 \cite{Kreyszig1991}. (At an inflection point the tangent line has a contact of order 2, in which case we say that the osculating circle has an infinity radius.) Further, we say that a curve $\alpha$ and a surface $\Sigma$ has a contact of order $k$ if there exists a curve in $\Sigma$ which has a contact of order $k$ with $\alpha$ and all the other curves has a lower, or equal, order of contact\footnote{For a level set surface $\Sigma=G^{-1}(c)$, an order $k$ contact is equivalent to $\beta^{(i)}(s_0^*)=0$ ($1\leq i\leq k$), where $\beta=G\circ\alpha$ and $c=\beta(s_0^*)=\alpha(s_0)$ \cite{Kreyszig1991}.}. For example, the osculating plane has a contact of order 2 with its reference curve, while the osculating sphere has a contact of order 3: when $\tau=0$, the osculating plane has a contact of order 3 and  the osculating sphere has an infinity radius. At a twisted point, the center and radius of the osculating sphere are respectively \cite{Kreyszig1991}
\begin{equation}
P_S = \alpha+\frac{1}{\kappa}\mathbf{n}+\frac{1}{\tau}\frac{{\rm d}}{{\rm d}s}\left(\frac{1}{\kappa}\right)\mathbf{b}\,\mbox{ and }\,R_S = \sqrt{\frac{1}{\kappa^2}+\frac{1}{\tau^2}\left[\frac{{\rm d}}{{\rm d}s}\left(\frac{1}{\kappa}\right)\right]^2}.
\end{equation}

\begin{Rem}
It is possible to prove Theorem \ref{thr::charSphCurves} using osculating spheres. This is useful when we do not have good orthogonality properties as, e.g., in spaces equipped with a degenerate metric \cite{daSilvaTJM2019}. However, the use of osculating spheres demands that the curve must be $C^4$ and also that $\tau\not=0$, while in Bishop's approach one needs just a $C^2$ condition and no restriction on the torsion. 
\end{Rem}

Using the concept of osculating spheres, we would intuitively say that every $C^4$ curve is locally spherical. In this case, it is tempting to ask whether the normals to the osculating spheres form an RM vector field. Unfortunately, this strategy does not work unless the curve is spherical. Indeed, using that $R_S'=\rho'\varsigma/\tau R_S$ and $P_S'=\varsigma\mathbf{b}$, where
\begin{equation}
\varsigma(s)=\tau(s)\rho(s)+\frac{{\rm d}}{{\rm d}s}\left(\frac{\rho'(s)}{\tau(s)}\right)\,\mbox{ and } \rho=\frac{1}{\kappa},\label{eq::SphAnnihilator}
\end{equation}
direct computation of the derivative of $\mathbf{N}=(\alpha-P_S)/R_S$ leads to
\begin{equation}
\frac{{\rm d}}{{\rm d}s}\left(\frac{\alpha(s)-P_S(s)}{R_S(s)}\right)=\frac{1}{R_S}\left[\mathbf{t}+\frac{\rho\rho'}{\tau R_S^2}\varsigma\,\mathbf{n}+\left(\frac{\rho'\,^2}{\tau^2R_S^2}-1\right)\varsigma\,\mathbf{b}\right].\label{eq::DerOfNormalToOscSphere}
\end{equation}
Then, the vector field  $\mathbf{N}=(\alpha-P_S)/R_S$ minimizes rotation if and only if $\alpha$ is spherical, i.e., when $\varsigma\equiv0$. (The condition to be spherical corresponds to $\varsigma\equiv0$ \cite{Kreyszig1991,Kuhnel2010}, which by direct examination of Eq. (\ref{eq::DerOfNormalToOscSphere}) is a necessary and sufficient condition to have $\mathbf{N}$ and $\mathbf{t}$ parallel.)

Now we investigate the possibility of extending $\kappa = \frac{1}{r}\sqrt{1+J^2}$ and $\tau = J'/(1+J^2)$, valid for spherical curves \cite{Kuhnel2010}, for a generic curve. Let $\Sigma_s = \mathbb{S}^2(P_S,R_S)$ be the osculating sphere of $\alpha$ at $\alpha(s)$. Near a fixed point $\alpha(s_0)$, we can obtain a spherical curve $\beta: (s_0-\epsilon,s_0+\epsilon)\to\Sigma_{s_0}$ by projecting $\alpha$ on $\Sigma_{s_0}$,
\begin{equation}
\beta(t)=r_0\frac{\alpha(t)-a_0}{\Vert\alpha(t)-a_0\Vert}\,,\label{eq::OscSphericalProjection}
\end{equation}
where $a_0=P_S(s_0)$ and $r_0=R_S(s_0)$. Since the osculating sphere has a contact of order 3 with the curve and $\kappa$ and $\tau$ are 2nd and 3rd order differential invariants, respectively, the torsion $\tau_{\alpha}$ and the curvature $\kappa_{\alpha}$ of a $C^3$ regular twisted curve $\alpha$ and the torsion $\tau_{\beta}$ and the curvature $\kappa_{\beta}$ of its (osculating) spherical projection $\beta$ coincide at $s_0$: 
$
\kappa_{\alpha}(s_0) = \kappa_{\beta}(s_0)$  and $\tau_{\alpha}(s_0) = \tau_{\beta}(s_0)$. Thus, we can write
\begin{equation}
\kappa_{\alpha}(s_0) = \frac{1}{R_S(s_0)}\sqrt{1+J^2(s_0)}\,,
\end{equation}
and 
\begin{equation}
\tau_{\alpha}(s_0) = \frac{\langle\alpha(s_0)-P_S(s_0),\alpha'(s_0
)\times\alpha'''(s_0)\rangle}{1+J^2(s_0)}=\frac{J'(s_0)}{1+J^2(s_0)}+\frac{\kappa(s_0)\varsigma(s_0)}{1+J^2(s_0)}\,,
\end{equation}
where $J(s)=\langle\alpha(s)-P_S(s),\alpha'(s)\times\alpha''(s)\rangle$ and $\varsigma$ is defined in Eq. (\ref{eq::SphAnnihilator}).

\begin{Rem}
From the expressions above we see that $\theta(s)=\arctan\,J(s)$ is only valid for spherical curves. In analogy with the study of the normals to the osculating spheres, the discrepancy between the results for a spherical and a generic curve is proportional to $\varsigma$, which vanishes only for a curve on a sphere.
\end{Rem}

\section{Rotation minimizing frames and spherical and plane curves}

Now, we address the problem of characterizing those curves whose position vector lies, up to a translation, on a moving plane spanned by their unit tangent and an RM vector field, i.e., curves such that
\begin{equation}
\alpha(s) - p = A(s)\,\mathbf{t}(s)+B(s)\,\mathbf{n}_1(s), 
\end{equation}
for some fixed point $p$ and some functions $A$ and $B$. This problem has to do with the more general quest of studying curves that lie on a given (moving) plane generated by two chosen vectors of a moving trihedron, such as normal and rectifying curves, i.e., curves that lie on their normal or rectifying planes, respectively. It is known that normal and rectifying curves correspond to spherical curves and geodesics on a cone \cite{ChenMonthly2003,ChenTJM2017}, respectively. Now, we have

\begin{thm}\label{theoCharPlnSphrCrvUsngMvngRMpln}
Up to a translation, the position vector of a $C^2$ regular curve $\alpha:I\to \mathbb{E}^3$ lies on a moving plane spanned by its unit tangent and a rotation minimizing vector field if and only if $\alpha$ is either a plane or a spherical curve.
\end{thm}
{\it Proof. } If $\alpha$ lies on an RM moving plane $\mbox{span}\{\mathbf{t},\mathbf{n}_1\}$, $
\alpha-p = A\mathbf{t}+B\mathbf{n}_1$, then
\begin{equation}
\mathbf{t}  =  (A'-\kappa_1B)\mathbf{t}+(B'+\kappa_1A)\mathbf{n}_1+\kappa_2A\,\mathbf{n}_2
      \Rightarrow \left\{
      \begin{array}{c}
      A'-\kappa_1B = 1\\
      B'+\kappa_1A = 0\\
      \kappa_2A    = 0\\
      \end{array}
      \right..\label{eq::coefRMplaneCurve}
\end{equation}
If $\kappa_2(s)=0$ for all $s$, then $\mathbf{n}_2$ is a constant vector and, consequently, $\alpha$ lies in the plane normal to $\mathbf{n}_2$. On the other hand, if $A(s)=0$ for all $s$, it follows from the second equation of (\ref{eq::coefRMplaneCurve}) that $B$ is a constant. In this case, $\alpha-p = B\,\mathbf{n}_1$ and $\alpha$ is spherical. (It lies on a sphere of radius $\vert B\vert$ and center $p$.)

Conversely, if $\alpha$ is a plane curve, then the normal $\mathbf{v}$ to the plane is constant and, consequently, should be RM: $\mathbf{v}'=0$. On the other hand, if the curve is spherical, the normal $\mathbf{v}=\frac{1}{R}(\alpha-p)$ minimizes rotation and, trivially,  $\alpha-p\in\mbox{span}\{\mathbf{t},\mathbf{v}\}$, with $\mathbf{v}$ an RM vector field.
\qed

\section{{Characterization of Bertrand curves and helices}}

The \emph{natural mate} of a regular  curve $\alpha:I\to\mathbb{E}^3$ is the curve $\beta$ such that $\mathbf{t}_{\beta}=\mathbf{n}$, i.e., $\beta(s)=\int\mathbf{n}(u){\rm d}u$. In Ref. \cite{DeshmukhTJM2018}, it is proposed the question ``Which Frenet curve has a spherical natural mate?". In this section, we  characterize curves with spherical natural mates in terms of Bertrand curves: $\alpha$ is said to be a \emph{Bertrand curve} if there exists a curve $\gamma$ such that $\alpha$ and $\gamma$ have the same principal normal. Bertrand curves are characterized by a linear relation between their curvature and torsion, see Theorem 23.1 of Ref. \cite{Kreyszig1991}. {In addition, we use the concept of mate curves to characterize \emph{slant helices}, i.e., curves whose principal normal makes a constant angle with a fixed direction \cite{IzumiyaTJM2004}, in terms of general helices, i.e., curves that make a constant angle with a fixed direction \cite{Kreyszig1991}.}

If $\{\mathbf{t},\mathbf{n},\mathbf{b}\}$ is the Frenet frame of $\alpha$, then we may consider $\{\mathbf{n},\mathbf{b},\mathbf{t}\}$ as an orthonormal moving frame along the natural mate $\beta$, whose equation of motion is
\begin{equation}\label{eq::EqsMotionNaturalMate}
    \frac{\rm d}{{\rm d} s}\left(
    \begin{array}{c}
         \mathbf{n} \\
         \mathbf{b} \\
         \mathbf{t} \\
    \end{array}
    \right)=
    \left(
    \begin{array}{ccc}
         0 & \tau & -\kappa \\
         -\tau & 0 & 0 \\
         \kappa & 0 & 0 \\
    \end{array}
    \right)
    \left(
    \begin{array}{c}
         \mathbf{n} \\
         \mathbf{b} \\
         \mathbf{t} \\
    \end{array}
    \right).
\end{equation}
Consequently, the Frenet frame leads to an RM frame for its natural mate. Thus, if $\kappa_{\beta}$ and $\tau_{\beta}$ denote the curvature and torsion of the natural mate, we can write
\begin{equation}
        \tau = \kappa_{\beta}\cos\theta,\,
        \kappa = -\kappa_{\beta}\sin\theta,\,\theta' = \tau_{\beta}.
\end{equation}
From these correspondences, we can devise an alternative proof for the known expressions for the curvature and torsion of a natural mate \cite{DeshmukhTJM2018}:
\begin{equation}
    \kappa_{\beta} = \sqrt{\kappa^2+\tau^2}\mbox{ and }\tau_{\beta}=\frac{\kappa^2}{\kappa^2+\tau^2}\left(\frac{\tau}{\kappa}\right)'.
\end{equation}
As observed by Deshmukh \textit{et al.} \cite{DeshmukhTJM2018}, it is straightforward to see that the natural mate is a plane curve if and only if $\alpha$ is a generalized helix, i.e., if and only if the ratio $\tau:\kappa$ is constant \cite{Kreyszig1991}: $\tau_{\beta}\equiv0\Leftrightarrow (\tau/\kappa)'\equiv0$. On the other hand, from Theorem \ref{thr::charSphCurves}, $\beta$ is spherical if and only if $\kappa$ and $\tau$ are linearly related, i.e.,  if and only if $\alpha$ is a Bertrand curve. Thus, we have the important characterization

\begin{thm}\label{thrCharBertrandCrv}
A regular curve $\alpha:I\to\mathbb{E}^3$ has a spherical natural mate if and only if it is a Bertrand curve. 
\end{thm}

Deshmukh \textit{et al.} showed that if $\alpha$ has constant curvature, then its natural mate is spherical. The converse was shown to be true if the natural mate is spherical but not planar, see Theorem 2 of Ref. \cite{DeshmukhTJM2018}. Taking into account Theorem \ref{thrCharBertrandCrv}, it is worth mentioning that curves with constant curvature have a Bertrand partner, see Theorem 23.3 of \cite{Kreyszig1991}.

{A curve $\alpha:I\to\mathbb{E}^3$ is a slant helix if and only if its principal normal $\mathbf{n}=\mathbf{t}_{\gamma}$ makes a constant angle with a constant direction $\mathbf{c}$. Therefore, we concluded that }

\begin{thm}\label{thrCharSlantHelix}
{
A regular curve $\alpha:I\to\mathbb{E}^3$ is a slant helix if and only if its natural mate is a general helix.
}
\end{thm}
{Moreover, using the characterization of general helices in terms of their torsion and curvature, we can recover the known characterization of slant helices in terms of their torsion and curvature \cite{IzumiyaTJM2004}:}
\begin{equation}
    { \frac{\kappa^2}{(\kappa^2+\tau^2)^{\frac{3}{2}}}\left(\frac{\tau}{\kappa}\right)'=\frac{\kappa_{\gamma}}{\kappa_{\gamma}}=\mbox{constant}. }
\end{equation}

We can extend the notion of mate curves to RM frames. More precisely, we say that $\gamma$ is an \emph{RM mate} of $\alpha$ if $\mathbf{t}_{\gamma}$ is an RM vector field along $\alpha$. Thus, equipping $\alpha$ with an RM frame,
\begin{equation}\label{eqMotionEqsRMmates}
   \frac{\mathrm{d}}{\mathrm{d}s}\left(
    \begin{array}{c}
         \mathbf{t} \\
         \mathbf{n}_{1} \\
         \mathbf{n}_{2} \\
    \end{array}
    \right)=
    \left(
    \begin{array}{ccc}
         0 & \kappa_{1} & \kappa_{2} \\
         -\kappa_{1} & 0 & 0 \\
         -\kappa_{2} & 0 & 0 \\
    \end{array}
    \right)
    \left(
    \begin{array}{c}
         \mathbf{t} \\
         \mathbf{n}_{1} \\
         \mathbf{n}_{2} \\
    \end{array}
    \right),
\end{equation}
 leads to a Frenet frame for its RM mate $\gamma$,
\begin{equation}
        \frac{\mathrm{d}}{\mathrm{d}s}\left(
    \begin{array}{c}
         \mathbf{n}_{1} \\
         -\mathbf{t} \\
         \mathbf{n}_{2} \\
    \end{array}
    \right)=
    \left(
    \begin{array}{ccc}
         0 & \kappa_{1} & 0 \\
         -\kappa_{1} & 0 & -\kappa_{2} \\
         0 & \kappa_{2} & 0 \\
    \end{array}
    \right)
    \left(
    \begin{array}{c}
         \mathbf{n}_{1} \\
         -\mathbf{t} \\
         \mathbf{n}_{2} \\
    \end{array}
    \right).
\end{equation}
Thus, the curvature and torsion of the RM mate are respectively given by $\kappa_{\gamma}=\kappa\cos(\int\tau)$ and $\tau_{\gamma}=-\kappa\sin(\int\tau)$. 

In a sense, RM mate is a dual concept of natural mate and, consequently, we should expect characterizations for RM mates dual to the characterizations of natural mates. Indeed, it follows from Eq. (\ref{eqMotionEqsRMmates}) and Theorem \ref{thr::charSphCurves} that
\begin{thm}\label{thrCharSphPlnCrvUsngRMmates}
A regular curve is a spherical (or plane) curve if and only if its RM mate is a Bertrand curve (or generalized helix, respectively).
\end{thm}

\section{Spherical curves in Lorentz-Minkowski space}

In the Lorentz-Minkowski space $\mathbb{E}_1^3$, i.e., $\mathbb{R}^3$ equipped with the indefinite  metric $\langle x,y\rangle_1=x_1y_1+x_2y_2-x_3y_3$, we have three types of spheres: hyperbolic planes $\mathbb{H}^2_0(p,r)$, de Sitter spaces $\mathbb{S}_1^2(p,r)$, and lightlike cones $\mathcal{C}^2(p)$ \cite{daSilvaJG2017,LopesIEJG2014}. 

Let $\alpha:I\to \mathbb{E}_1^3$ be a regular curve parameterized by arc-length $s$ and consider that $\alpha$ is a non-null curve on $\mathbb{H}_0^2(r)$ or $\mathbb{S}_1^2(r)$ (without loss of generality, the sphere is assumed to be centered at the origin, i.e., $p=0$). We may equip $\alpha$ with its Saban frame $\{\mathbf{t},\mathbf{u},\mathbf{t}\times_1\mathbf{u}\}$ \cite{SabanRendLincei1958}, where $\mathbf{u}=\alpha/r$. Adopting the notation $\epsilon_t=\langle\mathbf{t},\mathbf{t}\rangle_1\in\{-1,1\}$ and $\epsilon_u=\langle\mathbf{u},\mathbf{u}\rangle_1\in\{-1,1\}$, we may write $\alpha''$ as
\begin{eqnarray}
\alpha'' & = & \epsilon_t\langle\mathbf{t},\alpha''\rangle_1\,\mathbf{t}+\epsilon_u\langle\mathbf{u},\alpha''\rangle_1\,\mathbf{u}-\epsilon_t\epsilon_u\langle\mathbf{t}\times_1\mathbf{u},\alpha''\rangle_1\,\mathbf{t}\times_1\mathbf{u}\nonumber\\
 & = & -\frac{\epsilon_t\epsilon_u}{r^2}\,\alpha+\frac{\epsilon_t\epsilon_u}{r^2}J_1\,\alpha'\times_1\alpha,\label{eq::DecompositionAccelVectorLMspherCurv}
\end{eqnarray}
where $J_1=\langle\alpha,\alpha'\times_1\alpha''\rangle_1$. Supposing that $\mathbf{t}'$ is non-null and using that in this case $\alpha''=\epsilon_n\kappa\,\mathbf{n}$, we have (also notice that $\langle\alpha,\alpha\rangle_1=\mbox{const.}\Rightarrow\langle\alpha,\alpha''\rangle_1=-\epsilon_t$)
\begin{equation}
\kappa^2(s) = \frac{\epsilon_u\epsilon_n}{r^2}\Big(1-\epsilon_t\,J_1^2(s)\Big).
\end{equation}
Observe that, as happens in Euclidean space, the spherical curvature $J_1$ is just a constant multiple of the geodesic curvature of $\alpha$ on $\mathbb{H}_0^2(r)$ or $\mathbb{S}_1^2(r)$.

Now, let us compute  the torsion in terms of $J_1$. From the expressions $\mathbf{n}=\epsilon_n\alpha''/\kappa$, $\mathbf{b}=\epsilon_n\alpha'\times_1\alpha''/\kappa$, and $\tau=\langle \mathbf{n}',\mathbf{b}\rangle_1$, we have
\begin{equation}
\tau=-\langle \mathbf{n},\mathbf{b}'\rangle_1=-\langle\frac{\epsilon_n}{\kappa}\alpha'',\frac{\epsilon_n}{\kappa}\alpha'\times_1\alpha'''\rangle_1=\frac{\epsilon_t\epsilon_u}{r^2\kappa^2}\langle\alpha-J_1\alpha'\times_1\alpha,\alpha'\times_1\alpha'''\rangle_1.
\end{equation}
Now, using the identity
$
\langle\mathbf{A}\times_1\mathbf{B},\mathbf{C}\times_1\mathbf{D}\rangle_1 = \langle\mathbf{A},\mathbf{D}\rangle_1\langle\mathbf{B},\mathbf{C}\rangle_1-\langle\mathbf{A},\mathbf{C}\rangle_1\langle\mathbf{B},\mathbf{D}\rangle_1
$ for $\mathbf{A}=\alpha'$, $\mathbf{B}=\alpha$, $\mathbf{C}=\alpha'$, and $\mathbf{D}=\alpha'''$, we find that $\langle\alpha'\times_1\alpha,\alpha'\times_1\alpha'''\rangle_1=0$. Finally, the torsion of a non-null curve with non-null principal normal is
\begin{equation}
\tau(s)  = \frac{\epsilon_t\epsilon_u}{\kappa^2(s)r^2}J_1'(s) =\frac{J_1'(s)}{1-\epsilon_t\,J_1^2(s)}.
\end{equation}

Integration of the above equation gives $\tau=(\mathrm{arccoth}\,J_1)'$ if $\alpha$ is spacelike ($\epsilon_t=+1$) and $\tau=(\mathrm{arctan}\,J_1)'$ if $\alpha$ is timelike ($\epsilon_t=-1$). This is compatible with the fact that the normal plane of a spacelike (timelike) curve has to be timelike (spacelike) and, therefore, rotations are parameterized by hyperbolic (usual, respectively) trigonometric functions \cite{daSilvaJG2017,OzdemirMJMS2008}. On the other hand, if $\mathbf{t}=\alpha'$ is spacelike but $\mathbf{t}'=\alpha''$ is lightlike, then Eq. (\ref{eq::DecompositionAccelVectorLMspherCurv}) leads to $J_1=\pm r$ constant. It is worth mentioning that in such case $\alpha$ has to be a plane curve, see Remark 7 of \cite{daSilvaJG2017}. (This plane has to be lightlike, Proposition 4.2 of \cite{HondaArXiv2019}.)

\begin{Rem}
The characterization of plane and spherical curves in terms of the normal development curve $s\mapsto(\kappa_1(s),\kappa_2(s))$, see Theorem \ref{thr::charSphCurves}, is also valid in $\mathbb{E}_1^3$ \cite{daSilvaJG2017}. This is the key step to prove Theorems \ref{theoCharPlnSphrCrvUsngMvngRMpln}, \ref{thrCharBertrandCrv}, {\ref{thrCharSlantHelix}}, and \ref{thrCharSphPlnCrvUsngRMmates}. Thus, these theorems are valid in $\mathbb{E}_1^3$ as well. (For the study of natural mates in $\mathbb{E}_1^3$ see, e.g., \cite{ChoiJMAA2012}: notice Choi \emph{et al.} name the natural mates as principal-donor curves.)
\end{Rem}

\ack  The author would like to thank useful discussions with Gilson S. Ferreira Jr. (Universidade Federal Rural de Pernambuco) and the financial support provided by the Mor\'a Miriam Rozen Gerber Fellowship for Brazilian Postdocs.

\end{document}